\documentclass[12pt]{article}
\usepackage[ansinew]{inputenc}
\usepackage{array}
\usepackage{color}
\usepackage{amsmath}
\usepackage{amsxtra}
\usepackage{amstext}
\usepackage{amssymb}
\usepackage{latexsym}
\begin{document}

\sloppy
\newcommand{\proof}{{\it Proof~}}
\newtheorem{thm}{Theorem}[section]
\newtheorem{cor}[thm]{Corollary}
\newtheorem{lem}[thm]{Lemma}
\newtheorem{prop}[thm]{Proposition}
\newtheorem{eg}[thm]{Example}
\newtheorem{defn}[thm]{Definition}

\newtheorem{rem}[thm]{Remark}
\numberwithin{equation}{section}

\thispagestyle{empty}
\begin{center}
{\bf\Large\textbf{$(1+2u)$-constacyclic codes over $\mathbb{Z}_4+u\mathbb{Z}_4$}}\\
\vspace{.3cm}
{Mohammad Ashraf and Ghulam Mohammad}\\
\vspace{.2cm}
 Department of Mathematics\\
    Aligarh Muslim University \\
        Aligarh -202002(India)\\
{\it E-mails} :  {mashraf80@hotmail.com; mohdghulam202@gmail.com}
\end{center}

\parindent=0mm
\begin{abstract}
\noindent Let $R=\mathbb{Z}_4+u\mathbb{Z}_4,$ where $\mathbb{Z}_4$ denotes the ring of integers modulo $4$ and $u^2=0$. In the present paper, we introduce a new Gray map from $R^n$ to $\mathbb{Z}_{4}^{2n}.$ We study $(1+2u)$-constacyclic codes over $R$ of odd lengths with the help of cyclic codes over $R$. It is proved that the Gray image of $(1+2u)$-constacyclic codes of length $n$ over $R$ are cyclic codes of length $2n$ over $\mathbb{Z}_4$. Further, a number of linear codes over $\mathbb{Z}_4$ as the images of $(1+2u)$-constacyclic codes over $R$ are obtained.
\end{abstract}

\vspace{.4cm}
\parindent=7mm

\section{ Introduction} The study of codes over finite rings was initiated by Blake \cite{4, 5} in the year 1970. During the last decades of the twentieth century a great deal of attention has been given to codes over finite rings because of their new role in algebraic coding theory and their successful applications. In a landmark paper \cite{13}, it has been shown that certain good non-linear binary codes can be constructed from cyclic codes over $Z_4$ via the Gray map. Since then, codes over finite rings have been studied by many authors \cite{9, 10, 14, 17}. A lot of work has been done in this direction, but codes over $\mathbb{Z}_4$ remain a special topic of interest in the field of algebraic coding theory because of their relation to lattices, designs, cryptography and their many applications. There is a vast literature on codes over $\mathbb{Z}_4$ and their applications for detail see the references \cite{1, 7, 8, 12, 20, 21}.\\

Recently, Yildiz and Karadeniz considered linear and cyclic codes over the non-ring $F_2+uF_2+vF_2+uvF_2$ of size $16$ with $u^2=v^2=0$ and $uv=vu$ in \cite{22} and \cite{23}, where some good binary codes have been obtained as the images under two Gray maps. Motivated by this work the same authors \cite{24}, studied linear codes over the ring $\mathbb{Z}_4+u\mathbb{Z}_4,~u^2=0,$ which is also a non-chain ring of size $16.$ Further, Yildiz and Aydin \cite{25}, determined algebraic structure of cyclic codes over this ring. They obtained many new linear codes over $\mathbb{Z}_4$ as the Gray images of cyclic codes over $\mathbb{Z}_4+u\mathbb{Z}_4.$ Bandi and Bhaintwal \cite{6} gave the structural properties of cyclic codes over $\mathbb{Z}_4+u\mathbb{Z}_4$ of odd length in different aspect. They provided the general form of the generators of a cyclic code over $\mathbb{Z}_4+u\mathbb{Z}_4.$\\

Constacyclic codes constitute a remarkable generalization of cyclic codes, hence form an important class of linear codes in the coding theory. Constacyclic codes can be efficiently encoded using shift registers, which explains their preferred role in engineering. Qian et al. \cite{19}, introduced $(1+u)$-constacyclic codes of odd length over the ring $F_2+uF_2,$ where $u^2=0.$ with the help of cyclic codes over this ring. Later on Abualrub and Siap \cite{3}, studied $(1+u)$-constacyclic codes of arbitrary length over this ring and they proved that the Gray image of a $(1+u)$-constacyclic code is a binary cyclic code of length $2n.$ $(1+v)$-constacyclic codes of odd length over the ring $F_2+uF_2+vF_2+uvF_2$ were studied by Karadeniz and Yildiz in \cite{15}. A lot of work has been done on constacyclic codes over different structure of rings by several authors we refer to \cite{2, 11, 16, 18, 26} and the references therein. This is the motivation to study $(1+2u)$-constacyclic codes of odd length over $\mathbb{Z}_4+u\mathbb{Z}_4$ with $u^2=0.$\\

\vspace{.4cm}
\parindent=7mm

\section{Preliminaries} Let $\mathbb{Z}_4$ be the ring of integers modulo $4$. Consider the ring $R=\mathbb{Z}_4+u\mathbb{Z}_4=\{a+ub|~a,~b\in \mathbb{Z}_4\}$ with $u^2=0$. $R$ is a commutative ring with characteristic $4$ and it can be viewed as the quotient ring $\mathbb{Z}_4[u]/\langle u^2\rangle$. The units and non-units in $R$ can be characterized by $``a+bu$ is a unit in $R$ if and only if $a$ is a unit in $\mathbb{Z}_4$". Therefore $\{1, 3, 1+u, 1+2u, 1+3u, 3+u, 3+2u, 3+3u\}$ is a set of units of $R$ while $\{0, 2, u, 2u, 2+u, 2+2u, 3u, 2+3u\}$ is a set of non-units of $R$.\\

$R$ has a total of six ideals given by $\langle0\rangle=\{0\},~\langle 2u\rangle=\{0, 2u\},~\langle u\rangle=\{0, u, 2u, 3u\},~\langle 2\rangle=\{0, 2, 2u, 2+2u\},~\langle 2+u\rangle=\{0, 2u, 2+u, 2+3u\}$ and $\langle 2, u\rangle=\{0, 2, u, 2u, 3u, 2+u, 2+2u, 2+3u\}$. It is clear that $R$ is a local Frobenius ring with $\langle 2, u\rangle$ as its maximal ideal. The residue field is given by $R/\langle 2, u\rangle,$ which is isomorphic to binary field $F_2$.\\

A commutative ring $R$ is called a chain ring if its ideals form a chain under the relation of inclusion. From the ideals of $R$, we can see that they do not form a chain; for instance, the ideals $\langle u\rangle$ and $\langle 2\rangle$ are not comparable. Therefore, $R$ is a non-chain extension of $\mathbb{Z}_4$. Also $R$ is not a principle ideal ring; for example, the ideal $\langle 2, u\rangle$ is not generated by any single element of $R$.\\

Let $R^n$ be the set of all $n$-tuples over $R$, then a nonempty subset $C$ of $R^n$ is called a code of length $n$ over $R$. $C$ is called linear code of length $n$ over $R$ if it is an $R$-submodule of $R^n$. Each codeword $\overline c$ in such a code $C$ is just an $n$-tuple of the form $\overline{c}=(c_0, c_1, \cdots, c_{n-1})\in R^n $ and can be represented by a polynomial in $R[x]$ as follows:$$\overline{c}=(c_0, c_1, \cdots, c_{n-1})~\mbox{if and only if}~c(x)=\sum\limits_{i=0}^{n-1}c_ix^i\in R[x].$$

\vspace{.4cm}
\parindent=7mm

\section{Gray map over $R$} Here, we give the definition of the Gray map on $R^n$. Observe that any element $c\in R$ can be expressed as $c=a+ub$, where $a,~b\in \mathbb{Z}_4$. The Gray map $\phi:R\longrightarrow \mathbb{Z}_{4}^{2}$ is given by $\phi(c)=\phi(a+ub)=(b, 2a+b)$. This map can be extended to $R^n$ in a natural way as follows: $$\phi:R^n\longrightarrow \mathbb{Z}_{4}^{2n}$$
$$(c_0, c_1, ..., c_{n-1})\longmapsto(b_0, b_1, ..., b_{n-1}, 2a_0+b_0, 2a_1+b_1, ..., 2a_{n-1}+b_{n-1}),$$ where $c_i=a_i+ub_i,~0\leq i\leq n-1$.\\

\noindent Also the Lee weight on $R$ is defined as $$w_L(a+ub)=w_L(b, 2a+b),$$ where $w_L(b, 2a+b)$ denotes the usual Lee weight on $\mathbb{Z}_{4}^{2}.$ This weight can be extended to $R^n$ componentwise. The Lee weight of a codeword $\overline{c}=(c_0, c_1, ..., c_{n-1})\in R^n$ is the rational sum of the Lee weights of its components, that is, $w_L=\sum\limits_{i=0}^{n-1}{w_L(c_i)}$. For any $\overline{c}_1,~\overline{c}_2\in R^n$, the Lee distance $d_L$ is given by $d_L(\overline{c}_1, \overline{c}_2)=w_L(\overline{c}_1-\overline{c}_2)$. The minimum Lee distance of $C$ is the smallest nonzero Lie distance between all pairs of distinct codewords of $C$. The minimum Lee weight of $C$ is the smallest nonzero Lee weight among all codewords of $C$.\\

\noindent Now with the definitions of Gray map and Lee weight, we have the following obvious theorem:

\begin{thm} The map $\phi:R^n\longrightarrow \mathbb{Z}_{4}^{2n}$ is a distance preserving linear isometry. Thus, if $C$ is a linear code over $R$ of length $n$, then $\phi(C)$ is a linear code over $\mathbb{Z}_4$ of length $2n$ and the two codes have the same Lee weight enumerators.
\end{thm}

\vspace{.4cm}
\parindent=7mm

\section{$(1+2u)$-constacyclic codes of odd length over $R$} A cyclic shift on $R^n$ is a permutation $\sigma$ such that $$\sigma(c_0, c_1, \cdots, c_{n-1})=(c_{n-1}, c_0, c_1, \cdots, c_{n-2}).$$ A linear code $C$ over $R$ of length $n$ is called cyclic code if it is invariant under the cyclic shift $\sigma$, that is, $\sigma(C)=C.$\\

\noindent A $(1+2u)$-constacyclic shift $\tau$ on $R^n$ acts as $$\tau(c_0, c_1, \cdots, c_{n-1})=((1+2u)c_{n-1}, c_0, c_1, \cdots, c_{n-2}).$$ A linear code $C$ over $R$ of length $n$ is called $(1+2u)$-constacyclic code if it is invariant under the $(1+2u)$-constacyclic shift $\tau$, that is, $\tau(C)=C.$\\

Using the polynomial representation of codewords of $R^n$ in $R[x]$, we see that for a codeword $\overline c\in R^n,~\sigma(\overline{c})$ corresponds to $xc(x)$ in $R[x]/\langle x^n-1\rangle$ while $\tau(\overline{c})$ corresponds to $xc(x)$ in $R[x]/\langle x^n-(1+2u)\rangle$. The following propositions are the analogy of a well-known result for cyclic codes over finite fields. The proofs are also similar, therefore, we are omitting the proofs.

\begin{prop} A subset $C$ of $R^n$ is a linear cyclic code of length $n$ over $R$ if and only if its polynomial representation is an ideal of $R[x]/\langle x^n-1\rangle.$
\end{prop}

\begin{prop} A subset $C$ of $R^n$ is a linear $(1+2u)$-constacyclic code of length $n$ over $R$ if and only if its polynomial representation is an ideal of $R[x]/\langle x^n-(1+2u)\rangle.$
\end{prop}

It is noted that $(1+2u)^n=1+2u$ if $n$ is odd and $(1+2u)^n=1$ if $n$ is even. Therefore, we only study the properties of $(1+2u)$-constacyclic codes of odd length over $R$. Cyclic codes over $R$ of odd length are classified in \cite{25}. Using this classification we study $(1+2u)$-constacyclic codes over $R$ of odd length by introducing the following isomorphism from $R[x]/\langle x^n-1\rangle$ to $R[x]/\langle x^n-(1+2u)\rangle$:

\begin{thm} Let $\mu$ be the map of $R[x]/\langle x^n-1\rangle$ into $R[x]/\langle x^n-(1+2u)\rangle$ defined by $\mu(c(x))=c((1+2u)x)$. If $n$ is odd, then $\mu$ is a ring isomorphism.
\end{thm}
\noindent{\bf{Proof.}} Since $(1+2u)$ is a unit in $R$, $(1+2u)^n=1$. Also we know that if $n$ is odd, then $(1+2u)^n=1+2u.$ Now suppose $a(x)\equiv b(x)(\mbox{mod}x^n-1)$, that is, $a(x)-b(x)=(x^n-1)q(x)$ for some $q(x)\in R[x]$. Then
$$\begin{array}{lll}
a((1+2u)x)-b((1+2u)x)& = & ((1+2u)^nx^n-1)q((1+2u)x)\\
                     & = & ((1+2u)x^n-(1+2u)^2)q((1+2u)x)\\
                     & = & (1+2u)(x^n-(1+2u))q((1+2u)x),
\end{array}$$ which means if  $a(x)\equiv b(x)(\mbox{mod}x^n-1),$ then  $a((1+2u)x)\equiv b((1+2u)x)(\mbox{mod}x^n-(1+2u)).$ But the converse can easily be shown as well which means $$a(x)\equiv b(x)(\mbox{mod}x^n-1)\Leftrightarrow a((1+2u)x)\equiv b((1+2u)x)(\mbox{mod}x^n-(1+2u)).$$ Note that one side of the implication tells us that $\mu$ is well defined and the other side tells us that it is injective, but since the rings are finite this proves that $\mu$ is an isomorphism.\\

\noindent The following corollary is an immediate consequence of the above theorem:

\begin{cor} $I$ is an ideal of $R[x]/\langle x^n-1\rangle$ if and only if $\mu(I)$ is an ideal of\linebreak $R[x]/\langle x^n-(1+2u)\rangle$, where $n$ is odd.
\end{cor}

\noindent Before stating our next result, we need the following known lemma:

\begin{lem}{\cite[Theorem 4]{24}},  Let $n$ be odd and $C$ be a cyclic code of length $n$ over the ring $R$. Then $C$ is an ideal in $R[x]/\langle x^n-1\rangle$ generated by $$C=\langle a_1(x)(b_1(x)+2), ua_2(x)(b_2(x)+2)\rangle,$$ for some $a_i(x), b_i(x)\in \mathbb{Z}_4[x]$ such that $x^n-1=a_i(x)b_i(x)c_i(x)$ and $a_i(x), b_i(x), c_i(x)$ are monic coprime polynomials.
\end{lem}

\noindent Using the isomorphism $\mu$ and the above lemma, we characterize $(1+2u)$-constacyclic codes over $R$ of odd length as follows:

\begin{thm} Let $n$ be odd and $C$ be a $(1+2u)$-constacyclic code of length $n$ over the ring $R$. Then $C$ is an ideal in $R[x]/\langle x^n-(1+2u)\rangle$ generated by $$C=\langle a_1(\tilde{x})(b_1(\tilde{x})+2), ua_2(\tilde{x})(b_2(\tilde{x})+2)\rangle,$$ where $\tilde{x}=(1+2u)x$, $a_i(x), b_i(x)$ are the polynomials in $\mathbb{Z}_4[x]$ such that $x^n-1=a_i(x)b_i(x)c_i(x)$ and $a_i(x), b_i(x), c_i(x)$ are monic coprime polynomials.
\end{thm}

Now, we define a map ${\overline\mu}:R^n\longrightarrow R^n$ such that $${\overline\mu}(c_0, c_1, \cdots, c_{n-1})=(c_0, (1+2u)c_1, (1+2u)^2c_2, \cdots, (1+2u)^{n-1}c_{n-1}).$$ It is worth mentioning that ${\overline\mu}$ acts as the vector equivalent of $\mu$ on $R^n$. Therefore, we can restate Corollary 4.4 in terms of vectors as well.

\begin{cor} $C$ is a linear cyclic code over $R$ of odd length $n$ if and only if ${\overline\mu}(C)$ is a linear $(1+2u)$-constacyclic code of length $n$ over $R.$
\end{cor}

Note that if $c=a+ub\in R$, then $(1+2u)c=a+u(2a+b).$ Thus $$\begin{array}{lll}
w_L(c)&=&w_L(b, 2a+b)\\
       &=&w_L((1+2u)c).
\end{array}$$

\noindent In view of the Corollary 4.7, we have the following result:

\begin{cor} $C$ is a cyclic code over $R$ of length $n$ with Lee distance $d_L$ if and only if ${\overline\mu}(C)$ is a $(1+2u)$-constacyclic code over $R$ of length $n$ with same Lee distance, where $n$ is odd.
\end{cor}

\section{Gray images of $(1+2u)$-constacyclic codes over $R$} Even length cyclic codes over $\mathbb{Z}_4$ were characterized by Dougherty and Ling in \cite{12}. Here, we study even length cyclic codes over $\mathbb{Z}_4$ as the Gray images of $(1+2u)$-constacyclic codes over $R.$

\begin{prop}  Let $\tau$ be the $(1+2u)$-constacyclic shift of $R^n$ and $\sigma$ be the cyclic shift of $\mathbb{Z}_4^{2n}$. If $\phi$ is the Gray map from $R^n$ to $\mathbb{Z}_4^{2n},$ then $\phi\tau=\sigma\phi$.
\end{prop}
\noindent{\bf{Proof.}} Let $\overline{c}=(c_0, c_1, \cdots, c_{n-1})\in R^n$, where $c_i=a_i+ub_i$ with $a_i,~b_i\in \mathbb{Z}_4$ for $0\leq i\leq n-1$. Taking $(1+2u)$-constacyclic shift of $\overline{c}$, we have

$$\begin{array}{lll}
\tau(\overline{c}) & = & ((1+2u) c_{n-1}, c_0, \cdots, c_{n-2}) \\
      & = & ((1+2u)(a_{n-1}+ub_{n-1}), a_0+ub_0, \cdots, a_{n-2}+ub_{n-2}) \\
      & = & (a_{n-1}+u(2a_{n-1}+b_{n-1}), a_0+ub_0, \cdots, a_{n-2}+ub_{n-2}).
\end{array}$$ Now, using the definition of Gray map $\phi$, we can deduce that $$\phi(\tau(\overline{c}))=(2a_{n-1}+b_{n-1}, b_0, b_1, \cdots, b_{n-2}, b_{n-1}, 2a_0+b_0, 2a_1+b_1, \cdots, 2a_{n-2}+b_{n-2}).$$
On the other hand, $$\phi(\overline{c})=(b_0, b_1, \cdots, b_{n-1}, 2a_0+b_0, 2a_1+b_1, \cdots, 2a_{n-1}+b_{n-1}).$$
Hence, $$\sigma(\phi(\overline{c}))=(2a_{n-1}+b_{n-1}, b_0, b_1, \cdots, b_{n-2}, b_{n-1}, 2a_0+b_0, 2a_1+b_1, \cdots, 2a_{n-2}+b_{n-2}).$$
Therefore, $$\phi\tau=\sigma\phi.$$

\noindent As a consequence of Proposition 5.1, we get the following main result:

\begin{thm} The Gray image of a linear $(1+2u)$-constacyclic code over $R$ of length $n$ is a distance invariant linear cyclic code over $\mathbb{Z}_4$ of length $2n.$
\end{thm}
\noindent{\bf{Proof.}} Let $C$ be a linear $(1+2u)$-constacyclic code over $R$ of length $n.$ Then $\tau(C)=C,$ and therefore $\phi(\tau(C))=(\phi\tau)(C)=\phi(C).$ It follows from Proposition 5.1 that $(\sigma\phi)(C)=\sigma(\phi(C))=\phi(C),$ which means that $\phi(C)$ is a linear cyclic code over $\mathbb{Z}_4$ of length $2n.$\\

Yildiz and Aydin in \cite{25} obtained table for cyclic codes over $R$ of length $7.$ By modifying the generators with isomorphism $\mu$, we get the following table for $(1+2u)$-constacyclic codes of length $7$ over $R$. The generators of $(1+2u)$-constacyclic codes given by the Theorem 4.6 are taken to be $C=\langle g_1(\tilde{x}), ug_2(\tilde{x})\rangle$, where $\tilde{x}=(1+2u)x.$

\newpage

\noindent Some $(1+2u)$-constacyclic codes of length $7$ and their $\mathbb{Z}_4$ images.

\begin{center}
\begin{tabular}{|c|c|c|}
\hline
$g_1(\tilde{x})$ & $g_2(\tilde{x})$ & $\mathbb{Z}_4$ Parameters\\
\hline
$0$ & $3x^4+2x^3+x^2+3(1+2u)x+3$ & $[14, 4^32^0, 12]$ \\

$0$ & $x^4+(1+2u)x^3+3x^2+3$ & $[14, 4^32^3, 8]$ \\

$3x^4+2x^3+x^2+(3+2u)x+3$ & $x^4+(1+2u)x^3+3x^2+3$ & $[14, 4^62^3, 6]$ \\

$(3+2u)x^3+x^2+2x+1$ & $(1+2u)x^3+2x^2+(1+2u)x+1$ & $[14, 4^82^3, 4]$ \\

$(3+2u)x^3+x^2+2x+1$ & $(3+2u)x+1$ & $[14, 4^{10}2^0, 4]$ \\

$(3+2u)x^3+x^2+2x+1$ & $(1+2u)x+1$ & $[14, 4^{10}2^1, 4]$\\

$(1+2u)x^3+2x^2+(1+2u)x+1$ & $(1+2u)x+1$ & $[14, 4^{10}2^{4}, 2]$\\

$(3+2u)x+1$ & $(3+2u)x+1$ & $[14, 4^{12}2^0, 2]$\\

$(3+2u)x+1$ & $(1+2u)x+1$ & $[14, 4^{12}2^1, 2]$\\

$(1+2u)x+1$ & $(1+2u)x+1$ & $[14, 4^{12}2^2, 2]$\\

$(3+2u)x+1$ & $3$ & $[14, 4^{13}2^0, 2]$\\

\hline
\end{tabular}
\end{center}

\vspace{.6cm}

\noindent Now, we give some examples of $(1+2u)$-constacyclic codes over $R$ of odd lengths other than length $7$. Also, we find the $\mathbb{Z}_4$ images of these constacyclic codes.\\

\noindent\textbf{Example 5.3} Let $n=9$ and $g_1(x)=0,~g_2(x)=x^8+(1+2u)x^7+x^6+(1+2u)x^5+x^4+(1+2u)x^3+3x^2+(3+2u)x+3.$ Then $C=\langle g_1(x), ug_2(x)\rangle$ is a $(1+2u)$-constacyclic code of length $9$ over $R$ with minimum Lee distance $8.$ In view of Theorem 5.2, the Gray image $\phi(C)$ of $C$ is a cyclic code over $\mathbb{Z}_4$ with parameters $[18, 4^12^6, 8]$.\\

\noindent\textbf{Example 5.4} Let $n=15$ and $g_1(x)=2x^{10}+2x^{8}+2x^{5}+2x^{4}+2x^2+2x,~g_2(x)=(3+2u)x^{13}+x^{12}+3x^{10}+(1+2u)x^9+(3+2u)x^7+x^6+3x^4+(1+2u)x^3+(3+2u)x+1.$ Then $C=\langle g_1(x), ug_2(x)\rangle$ is a $(1+2u)$-constacyclic code of length $15$ over $R$ with minimum Lee distance $8.$ In view of Theorem 5.2, the Gray image $\phi(C)$ of $C$ is a cyclic code over $\mathbb{Z}_4$ with parameters $[30, 4^22^{10}, 8]$.\\

\noindent\textbf{Example 5.5} Let $n=23$ and $g_1(x)=0,~g_2(x)=x^{22}+(1+2u)x^{21}+x^{20}+(1+2u)x^{19}+x^{18}+(1+2u)x^{17}+x^{16}+(1+2u)x^{15}+x^{14}+(1+2u)x^{13}+x^{12}+(3+2u)x^{11}
+3x^{10}+(1+2u)x^{9}+x^8+(1+2u)x^7+3x^6+(3+2u)x^5+3x^4+(1+2u)x^3+3x^2+(1+2u)x+3.$ Then $C=\langle g_1(x), ug_2(x)\rangle$ is a $(1+2u)$-constacyclic code of length $23$ over $R$. In view of Theorem 5.2, the Gray image $\phi(C)$ of $C$ is a cyclic code over $\mathbb{Z}_4$ with parameters $[46, 4^12^{11}, 28]$.\\

One-generator cyclic codes over $R$ of length $n$ have been studied by Yildiz and Aydin in \cite{25}. Motivated by this study we consider a one generator $(1+2u)$-constacyclic code over $R$ of length $n$, that is, we take an ideal in $R[x]/\langle x^n-(1+2u)\rangle$ generated by some polynomial $a(x)+ub(x)\in R[x]/\langle x^n-(1+2u)\rangle,$ where $a(x),~b(x)\in \mathbb{Z}_4[x]$ and deg$(a(x))<n,$ deg$(b(x))<n.$ In light of this, we give the following  result about the Gray image of such a code in $\mathbb{Z}_4$:\\

\noindent\textbf{Theorem 5.6} Let $C=\langle a(x)+ub(x)\rangle$ be a $(1+2u)$-constacyclic code over $R$ of length $n.$ Then $\phi(C)$ is a cyclic code over $\mathbb{Z}_4$ of length $2n$ generated by the polynomial $b(x)+x^n(2a(x)+b(x))$ and $a(x)+x^na(x).$\\

\noindent{\bf{Proof.}} First we define a Gray map for polynomials as follows: $$\phi:R[x]/\langle x^n-(1+2u)\rangle\longrightarrow\mathbb{Z}_{4}[x]/\langle x^n-1\rangle\times\mathbb{Z}_{4}[x]/\langle x^n-1\rangle$$
$$\phi(a(x)+ub(x))=(b(x), 2a(x)+b(x)).$$ One can easily proved that $\phi$ is well defined. It may be noted that $(b(x), 2a(x)+b(x))$ gives us the same vector as $b(x)+x^n(2a(x)+b(x))$ in $\mathbb{Z}_{4}[x]/\langle x^{2n}-1\rangle.$ Since $C$ is an ideal in $\mathbb{Z}_{4}[x]/\langle x^{n}-(1+2u)\rangle,$ $ur(x)(a(x)+ub(x))\in C$ for all $r(x)\in \mathbb{Z}_4[x],$ and therefore we get $ur(x)a(x)\in C.$ Also $\phi(ur(x)a(x))=r(x)(a(x), a(x)),$ which gives us the same vector as $r(x)(a(x)+x^nr(x))$ in $\mathbb{Z}_{4}[x]/\langle x^{2n}-1\rangle.$ Again it may be noted that $$(r(x)+uq(x))(a(x)+ub(x))=r(x)(a(x)+ub(x))+uq(x)a(x).$$ Thus, we have
$$\phi[(r(x)+uq(x))(a(x)+ub(x))]=r(x)[b(x)+x^n(2a(x)+b(x))]+q(x)[a(x)+x^na(x)].$$ Since $r(x)$ and $q(x)$ are arbitrary polynomials in $\mathbb{Z}_4[x]$ of degree $<n,$ this proves the theorem.\\

\noindent We close our discussion with the following example of even length $(1+2u)$-constacyclic code over $R$:\\

\noindent\textbf{Example 5.7} Let $n=6$ and the one-generator $(1+2u)$-constacyclic code $C=\langle (x+x^2+x^4+x^5)+u(2x+2x^3+2x^5)\rangle.$ Then by Theorem 5.6, $\phi(C)$ is a cyclic code over $\mathbb{Z}_4$ of length $12$ with Lee distance $8$ generated by the polynomials $2x^{11}+2x^{10}+2x^9+2x^8+2x^5+2x^3+2x$ and $x^{11}+x^{10}+x^8+x^7+x^5+x^4+x^2+x.$\\

\section{Conclusion} In this paper, we have studied $(1+2u)$-constacyclic codes of odd length over $R=\mathbb{Z}_4+u\mathbb{Z}_4,$ where $u^2=0.$ We have defined a new Gray map from $R^n$ to $\mathbb{Z}_4^{2n}$ to study linear codes over $\mathbb{Z}_4$ as the Gray images of $(1+2u)$-constacyclic codes over $R.$ Further, we have proved that the Gray image of $(1+2u)$-constacyclic codes of length $n$ over $R$ are cyclic codes of length $2n$ over $\mathbb{Z}_4$.\\

We have used the fact that $(1+2u)^n=1$ when $n$ is even and $(1+2u)^n=(1+2u)$ when $n$ is odd, these two conditions are also true for the unit $(3+2u)$ in $R.$ Therefore, the study that we have done for $(1+2u)$-constacyclic codes of odd length over $R$ can be obtained analogously for $(3+2u).$ Also, one may generalized this study of constacyclic codes over the large ring $\mathbb{Z}_q+u\mathbb{Z}_q,$ where $q=p^m.$

\vspace{.3cm}

\begin{center}

\end{center}

\end{document}